\documentclass[a4paper,12pt,twoside]{amsart}
\usepackage{amsmath,amsthm,amsfonts}
\usepackage[T1]{fontenc}

\newtheorem{thm}{Theorem}[section]
\newtheorem{cor}[thm]{Corollary}
\newtheorem{lem}[thm]{Lemma}
\newtheorem{lemma}[thm]{Lemma}
\newtheorem{prop}[thm]{Proposition}
\newtheorem{proposition}[thm]{Proposition}
\newtheorem{definition}[thm]{Definition}

\theoremstyle{definition}

\newtheorem{rem}{Remark}

\newcommand{\Cal}{\mathcal}
\newcommand{\cal} {\mathcal}
\newcommand{\R}{{\mathbb{R}}}

\newcommand{\Q}{{\mathbb{Q}}}

\newcommand{\Z}{{\mathbb{Z}}}
\newcommand{\N}{{\mathbb{N}}}
\newcommand{\T}{{\mathbb{T}}}

\def\eop{\qed}
\def\Proof {\vskip -2mm {{\it Proof}.}}
\def\proof {\vskip -2mm {{\it Proof}.}}

\def \mod {\text{ mod }}
\def \ev {\cal E}
\newcommand{\moda} { {\rm \ mod \ } 1}

\title [Remarks on step cocycles over rotations]
{Remarks on step cocycles over rotations, \\ centralizers and
coboundaries}

\subjclass[2000] {28D05, 37A20, 37A40, 11Jxx} \keywords{cylinder
map, non regular cocycle, coboundary, diophantine properties}

\begin{document}

\baselineskip 15pt
\parindent=0mm

\maketitle \centerline {Jean-Pierre Conze and Jonathan Marco }
\vskip 3mm \centerline {IRMAR, University of Rennes 1}

\begin{abstract}
By using a cocycle generated by the step function $\varphi_{\beta,
\gamma} = 1_{[0, \beta]} - 1_{[0, \beta]} (. + \gamma)$ over an
irrational rotation $x \to x + \alpha \mod 1$, we present examples
which illustrate different aspects of the general theory of cylinder
maps. In particular, we construct non ergodic cocycles with ergodic
compact quotients, cocycles generating an extension $T_{\alpha,
\varphi}$ with a small centralizer. The constructions are related to
diophantine properties of $\alpha, \beta, \gamma$.
\end{abstract}

\tableofcontents

\vskip 3mm \section*{\bf Introduction}

Skew maps (also called cylindrical systems) yield an important
source of examples of dynamical systems preserving an infinite
invariant measure. In particular the class of skew maps over
1-dimensional irrational rotations using a step function as skewing
function has been widely studied in the literature. (cf.
\cite{Or83}, \cite{LePaVo96}, \cite{ALV98} and for other references
\cite{Co09}).

Our main examples here will be the cocycles generated over an
irrational rotation $T_\alpha: x \to x + \alpha \mod 1$ by the step
functions\footnote{In what follows the arguments of the functions
are taken modulo 1.}
$$\varphi_{\beta}(x) := 1_{[0,\beta]}
(x)-\beta, \ \varphi_{\beta, \gamma}(x) := 1_{[0,\beta]}
(x)-1_{[0,\beta]}(x+\gamma).$$

This simple function can be used to answer natural questions about
cocycles. In particular, we are interested in the construction of
non ergodic cocycles with ergodic compact quotients and cocycles
generating an extension $T_{\alpha, \varphi}: (x, y) \to (x+\alpha,
y + \varphi(x))$ with a small centralizer. This has the advantage to
illustrate the general ergodic theory of dynamical systems in
infinite measure through a very elementary and natural object.

After reminders on extensions of dynamical systems, essential values
and regularity of cocycles, we discuss some issues on
$\Z^2$-cocycles and centralizer of cylindrical maps. Then we present
general results on coboundaries equations over rotations and recall
results of M. Gu\'enais and F. Parreau on a multiplicative
quasi-coboundary equation. In the case of step functions, we give
sufficient conditions for solving in $L^2(\T^1)$ the linear
coboundary equation for the function $1_{[0, \beta]} - T_\gamma
1_{[0, \beta]}$.

As a result, it follows (Theorem \ref{ergoComp}) that there are real
numbers $\beta$ such that: \hfill \break - on one hand, for almost
every $\gamma$ the cocycle defined by $\varphi_{\beta,\gamma}$ is
non~regular (in particular it is not ergodic, but not a coboundary),
but all the compact quotients of the associated skew product are
ergodic, \hfill \break - on the other hand, there is an uncountable
set of values of $\gamma$ for which $\varphi_{\beta,\gamma}$ is a
coboundary.

Then we show different kinds of centralizer for $T_{\alpha,
\varphi_\beta}$: non trivial uncountable (case of unbounded partial
quotient), trivial (case of bounded quotients). At the opposite we
investigate also a property of "rigidity" for $\alpha$ of bounded
type, with an example of cocycle $\varphi$ which generates an
extension $T_{\alpha, \varphi}$ with a small centralizer. A last
application is the construction of a counter example in a conjugacy
problem for a group family. In the appendix, under diophantine
conditions on $\beta, \gamma$, we solve the linear coboundary
equation for $\varphi_{\beta, \gamma}$.

\vskip 3mm The authors are grateful to M. Lema\'nczyk for references
and comments on the centralizer, as well as to the referee for his
numerous and very helpful suggestions.

\vskip 5mm \newpage
\section{\bf Preliminaries on cocycles}
\subsection{Cocycles and group extension of dynamical systems}

\

In these preliminaries, we recall some standard facts on skew
products and regular cocycles.

Let $(X, {\Cal A}, \mu, T)$ be a dynamical system, i.e., a
probability space $(X, {\Cal A}, \mu)$ and a measurable invertible
transformation $T$ of $X$ which preserves $\mu$. In the sequel we
will {\it assume $T$ ergodic}.  Let $\varphi: X \to G$ be a
measurable function from $X$ to an abelian locally compact second
countable (lcsc) group $G$, with $m$ or $m_G$ denoting the Haar
measure on $G$.

The {\it skew product (or cylinder map)} over $(X, \mu, T)$ with the
fiber $G$ and the displacement (or skewing) function $\varphi$ is
the dynamical system $(X \times G, \mu \otimes m, T_\varphi)$, where
$$T_\varphi(x,g) = (T x, g + \varphi(x)).$$

For $n \in \Z$ we have $T_\varphi^n(x,g)=(T^n x, g+\varphi_n(x))$,
where $(\varphi_n)$ is the associated cocycle generated by $\varphi$
over the dynamical system:
\begin{equation} \label{cocycle_eq}
\varphi_n(x) = \sum_{j=0}^{n-1} \varphi(T^j x), \ n \geq 1.
\end{equation}

 For simplicity, the function $\varphi$ itself will be called a {\it
cocycle}. We say that a cocycle $\varphi : X \rightarrow G$ is
ergodic if the transformation $T_\varphi$ is ergodic on $X \times G$
for the measure $\mu \times m_G$.

\vskip 3mm Recall that two cocycles $\varphi$ and $\psi$ over a
dynamical system $(X, \mu, T)$ are {\it cohomologous} with {\it
transfer function} $\eta$, if there is a measurable map $\eta : X
\rightarrow G$ such that\footnote{If $f$ is a function defined on a
space $X$ and $T$ a transformation on $X$, we write simply $Tf$ for
the composed function $f\circ T$. The equalities between functions
are understood $\mu$-a.e.}
\begin{equation}
\varphi = \psi + T\eta - \eta. \label{cohom}
\end{equation}
$\varphi$ is a {\it $\mu$-coboundary} if it is cohomologous to 0.

\vskip 3mm {\em Recurrence}: When $G$ is non compact, to deal with
extensions with a non dissipative behavior, it is desirable that a
recurrence property holds. A point $x\in X$ is {\em recurrent} for
the cocycle $\varphi$, if $\varphi_n(x) \not \to \infty$ when $n$
tends to $\infty$. We say that $\varphi$ is {\em recurrent} if a.e.
$x \in X$ is recurrent. If the cocycle is recurrent, then the map
$T_\varphi$ is conservative for the invariant $\sigma$-finite
measure $\mu \times m_G$.

An integrable cocycle $\varphi$ with values in $\R$ is recurrent if
and only if $\int \varphi \ d\mu = 0$ (cf. \cite{Sc77}). If
$\varphi$ is a recurrent cocycle, than every cocycle cohomologous to
$\varphi$ is recurrent.

\vskip 3mm
\subsection{Essential values, non~regular cocycle}

\

First we recall the notion of essential values of a cocycle (cf. K.
Schmidt \cite{Sc77}, see also J. Aaronson \cite{Aa97}).

Let $\varphi$ be a cocycle with values in an abelian lcsc group $G$.
If $G$ is a non compact group, we add to $G$ a point at $\infty$
with the natural notion of neighborhood.

\begin{definition}\label{valess} {\rm An element $a \in G \cup \{\infty\}$
is an {\it essential value} of the cocycle $\varphi$ (over the
system $(X, \mu, T)$) if, for every neighborhood $V(a)$ of $a$, for
every measurable subset $B$ of positive measure,
\begin{eqnarray}
\mu(B\cap T^{-n} B \cap \{x: \varphi_n(x) \in V(a) \}\bigr) > 0,
{\rm \ for \ some \ } n \in \Z. \label{visitVa}
\end{eqnarray} } \end{definition}

We denote by ${\overline \ev}(\varphi)$ the set of essential values
of the cocycle $\varphi$ and by $\ev(\varphi) = {\overline
\ev}(\varphi)\cap G$ the {\it set of finite essential values}.

The set $\ev(\varphi)$ is a closed subgroup of $G$, with ${\cal
E}(\varphi) = G$ if and only if $(X \times G, \mu \otimes m,
T_\varphi)$ is ergodic.

Two cohomologous cocycles have the same set of essential values.
$\varphi$ is a coboundary if and only if ${\overline {\cal E}}
(\varphi) =\{ 0 \}$.

\begin{definition}\label{regul0} {\rm We say that the cocycle defined
by $\varphi$ is {\it regular}, if $\varphi$ can be reduced by
cohomology to an ergodic cocycle $\psi$ with values in the closed
subgroup $\ev(\varphi)$:}
\begin{eqnarray}
\psi = \varphi + \eta - T\eta, \label{regulDef}
\end{eqnarray}
\end{definition}

Let us recall some of the properties of regular cocycles. A cocycle
$\varphi$ is regular if and only if $\varphi / \ev(\varphi)$ is a
coboundary. A regular cocycle is recurrent. In the regular case
there is a "nice" ergodic decomposition of the measure $\mu \times
m_G$ for the skew map $T_\varphi$: any $T_\varphi$-invariant
function can be written $F(y - \eta(x))$ for a function $F$ which is
invariant by translation by elements of $\ev(\varphi)$, with $\eta$
given by (\ref{regulDef}).

If the cocycle is non~regular, then the measures $\mu_x$ on $X$ on
which is based the ergodic decomposition of $\mu \otimes m$ are
infinite, singular with respect to the measure $\mu$ and there are
uncountably many of them pairwise mutually singular (cf. K. Schmidt
\cite{Sc77}, see also \cite{CoRa09} for a complete description of
the ergodic decomposition in the general case of non abelian lcsc
groups $G$).

The following lemma is a simple tool which can be used to construct
non~regular cocycles.

\begin{lem}\label{sgroupdis} If $\varphi$ is a $\Z$-valued cocycle
such that there exists $s \not \in \Q$ for which the multiplicative
coboundary equation $e^{2\pi i s \varphi} = \psi / T\psi$ has a
measurable solution $\psi$, then ${\cal E}(\varphi) = \{0\}$. If
$\varphi$ is not a coboundary, then ${\overline {\cal E}}(\varphi) =
\{0, \infty \}$ and the cocycle $\varphi$ is non~regular.
\end{lem}
\proof \ From the hypothesis we have $\varphi = s^{-1} \zeta + \eta
- T\eta$, where $\zeta$ has values in $\Z$. The cocycle $\varphi$
can be viewed as a real cocycle with values in $\Z$, which is
cohomologous to a cocycle with values in the closed subgroup $s^{-1}
\Z$, with $s^{-1} \not \in \Q$.

In general, if a cocycle $\varphi$ is cohomologous to $\varphi_1$
and to $\varphi_2$, two functions with values respectively in closed
subgroups with an intersection reduced to $\{0\}$, then ${\cal
E}(\varphi) = {\cal E}(\varphi_1) \cap {\cal E}(\varphi_2) = \{0\}$.
\eop

\vskip 3mm

{\bf Cocycles and ergodicity in compact quotients}

If $G$ is compact, then there exist a measurable function $\eta : X
\to G$ such that, for the cocycle $\psi(x) = \varphi(x) + \eta(x) -
\eta(Tx) \in \ev(\varphi)$, the map $T_\psi$ is ergodic on $X \times
\ev(\varphi)$. Therefore $T_\varphi$ is regular.

Ergodicity implies ergodicity for all compact quotients $X \times G
/ G_0$, where $G_0$ is any cocompact closed subgroup of $G$. The
converse does not hold in general.

A question is to find examples of skew products which are non
ergodic on $X \times G$, but ergodic on all compact quotients $X
\times G/ G_0$.

There are example of skew products for which all compact quotients
are ergodic. For instance, the directional billiard in the plane
with periodic rectangular obstacles yields such examples: for almost
every direction the compact quotients of the directional billiard
are ergodic; nevertheless, due to recent results of K. Fr\k{a}czek
and C. Ulcigrai (\cite{FrUl11}), it is known that the billiard map
is non ergodic and even non~regular for a.e. parameters. This
provides examples, but we would like to construct more elementary
explicit examples (see Subsection \ref{nonErgCQt}).

\begin{rem} \label{rema1} Let $\varphi$ be a cocycle with values in
$G = \R^d \times \Z^{d'}$. If all of its compact quotients are
ergodic, then $\varphi$ is ergodic or non~regular. Indeed, if
$\varphi$ is regular, then $\varphi/ \ev(\varphi)$ is a coboundary.
Hence, if the compact quotients are ergodic for $\varphi$, the
compact quotients of $G / \ev(\varphi)$ are trivial. This implies
$\ev(\varphi) = G$ and $\varphi$ is ergodic.
\end{rem}

\vskip 3mm \subsection{$\Z^2$-actions and centralizer}

\

\vskip 3mm \subsubsection{$\Z^2$-actions and skew maps}
\label{Z2cocycle}

\

The construction of skew maps can be extended to group actions
generalizing the action of $\Z$ generated by iteration of a single
automorphism. We consider the case of $\Z^2$-actions.

Let $T_1, T_2$ be two commuting measure preserving invertible
transformations on $(X, \mu)$. They define a $\Z^2$-action on $(X,
\mu)$. A $G$-valued function $\varphi(n_1,n_2, x)$ on $\Z^2 \times
X$ is a cocycle for this action, if it satisfies the cocycle
relation:
$$\varphi(n_1+n_1',n_2+n_2', x) = \varphi(n_1,n_2, x) + \varphi(n_1',n_2',
T_1^{n_1} T_2^{n_2} x), \forall \, n_1, n_1',n_2, n_2' \in \Z.$$

Let $\varphi_i$, $i= 1,2$, be two measurable $G$-valued functions on
$X$ and consider the skew products $\tilde T_i: (x,y)
\longrightarrow (T_i x, y + \varphi_i\,(x))$ on $X \times \R$. Do
they generate a $\Z^2$-action which extends the $\Z^2$-action on
$(X, \mu)$?

The maps $\tilde T_1$ and $\tilde T_2$ commute if and only if the
following coboundary equation is satisfied
\begin{eqnarray}
\varphi_1 - T_2{\varphi_1} = \varphi_2 - T_1
{\varphi_2}. \label{commute}
\end{eqnarray}

If (\ref{commute}) is satisfied, then the composed transformation
$\tilde T_2^{n_2} \tilde T_1^{n_1}$ reads:
$$\tilde T_2^{n_2} \tilde T_1^{n_1}(x,y) = (T_2^{n_2}
T_1^{n_1}x, y + \varphi(n_1,n_2,x)),$$ where $\varphi(1, 0, x) =
\varphi_1(x)$, $\varphi(0, 1, x) = \varphi_2(x)$,
$\varphi(n_1,n_2,x)$ satisfies the cocycle relation and $(n_1, n_2)
\to \tilde T_2^{n_2} \tilde T_1^{n_1}$ defines a measure preserving
action of $\Z^2$ on $X \times G$.

Therefore it is equivalent to find $G$-valued $\Z^2$-cocycles (and
the corresponding skew products) or to find pairs $(\varphi_1,
\varphi_2)$ satisfying (\ref{commute}).

Clearly if $\varphi_1 = v - T_1\,v $ for some measurable function
$v$, then Equation (\ref{commute}) holds with $\varphi_2 = v -
T_2\,v$. A question is the construction of a pair $(\varphi_1,
\varphi_2)$ which satisfies (\ref{commute}) but are not of this
form. In other words, can we {\it construct solutions of
(\ref{commute}) which are not coboundaries?}

The answer depends on the choice of the transformations and on the
class to which the functions $\varphi_1, \varphi_2$ belong. For
instance, there is a "rigidity" for the $\Z^2$-shift on
$\{0,1\}^{\Z^2}$ endowed with the product measure. When the
functions are locally constant, the only solutions in that case are
the trivial ones (cf. K. Schmidt \cite{Sc95}, O. Jenkinson
\cite{Je01}).

In the case of rotations on $\mathbb{T}^1$, by using Fourier
analysis, we will give below explicit examples of non degenerate
solutions of (\ref{commute}) in $L^2(\T^1)$ (Theorem \ref{majBeta4})
and apply it to the construction of non trivial centralizers, a
notion that we recall now.

\vskip 3mm \subsubsection{Centralizer of the cylinder product}
\label{subsect-centralz}

\

A problem related to the construction of $\Z^2$-cocycles is the
study of the centralizer.

In what follows\footnote{The centralizer, in a wider sense, is the
collection of non-singular transformations of X which commute with
$\tilde T_1$ (see for instance \cite{ALV98}).} by {\it centralizer}
of a cylinder map $\tilde T_1: (x,y) \to (T_1 x, y + \varphi_1)$, we
mean the group $\Cal C(\tilde T_1)$ of {\it measure preserving
automorphisms} of $(X \times G, \mu \times m_G)$ which commute with
$\tilde T_1$. It contains the powers of the map and the translations
on the fibers. The skew products of the form $(x,y) \to (T_2 x, y +
\varphi_2)$ with $T_2$ commuting with $T_1$ and $(\varphi_1,
\varphi_2)$ satisfying (\ref{commute}) are elements of the group
$\Cal C(\tilde T_1)$.

\subsection{Case of an irrational rotation\label{1.4}}

\

In this subsection, we take the dynamical system $(X, \mu, T)$ in
the class of rotations on $\T^1$ (which could be replaced by a
compact abelian group $K$). For simplicity we take cocycles with
values in $\R$. About the centralizer and related questions, see
\cite{Kw92} (for $G$ a compact group), \cite{LeM97}, \cite{ALMN95},
\cite{ALV98}, \cite{LeMeNa03}.

In the sequel $\alpha$ will be an irrational number and $T_\alpha$
the ergodic rotation $x \to x+\alpha \text{ mod } 1$ on $X = \T^1$.
For a measurable function $\varphi : X \to \R$, we consider the skew
product $T_{\alpha, \varphi}: (x,y) \to (x + \alpha, y +
\varphi(x))$.

In this case, according (\ref{commute}), the automorphisms given by
skew products of the form $T_{\gamma, \psi}: (x,y) \to (x+\gamma, y
+ \psi)$, for $\gamma \in \T^1$ and a measurable function $\psi$
with $(\varphi, \psi)$ satisfying $\varphi - T_\gamma \varphi = \psi
- T_\alpha \psi$ are elements of the group $\Cal C(T_{\alpha,
\varphi})$. A problem is to find all elements in $\Cal C(T_{\alpha,
\varphi})$.

\vskip 3mm The following result is a special case of Proposition 1.1
in \cite{ALMN95}.

\begin{thm} \label{thmCommut} (cf. \cite{ALMN95}) Suppose that the cocycle
generated by $\varphi$ over the rotation $T_\alpha$ is ergodic. Then
any automorphism of $(X \times \R, \mu \times dy)$ commuting with
$T_{\alpha, \varphi}$ has the form $(x, y) \to (x+\gamma,
\varepsilon y +\psi(x))$ where $\gamma \in \T^1$, $\varepsilon$ is a
constant in $\pm 1$ and $\psi: X \to \R$ is a measurable function
such that
\begin{eqnarray}
\varepsilon \varphi - T_\gamma{\varphi} = \psi - T_\alpha \psi.
\label{commute2}
\end{eqnarray}
\end{thm}
\proof \  We give a sketch of the proof. The measure theoretic
details are omitted. Let $\tilde T_2$ be an automorphism which
commutes with $\tilde T_1 := T_{\alpha, \varphi}$. With the notation
$u(x,y) = e^{2\pi i x}$, we deduce from the commutation $\tilde T_1
\tilde T_2 = \tilde T_2 \tilde T_1$, that $\tilde T_2 u$ is an
eigenfunction for $\tilde T_1$ with eigenvalue $e^{2\pi i \alpha}$.
By ergodicity of $\tilde T_1$, this implies that $u \circ \tilde T_2
= \lambda u$, where $\lambda$ is a complex number of modulus 1.

It follows that $\tilde T_2$ leaves invariant the rotation factor of
$\tilde T_1$ and that there are $\gamma \in \R$ and a measurable map
$(x,y) \to V(x,y)$ from $X \times \R$ to $\R$ such that $\tilde T_2$
can be represented as $(x,y) \to \tilde T_2(x,y) = (x+\gamma,
V(x,y))$.

The commutation of the maps $\tilde T_1, \tilde T_2$ implies:
\begin{eqnarray}
V(x+\alpha, y + \varphi(x)) = V(x, y )+ \varphi(x+\gamma).
\label{Vphi}
\end{eqnarray}

Let us define $u_z(x,y):= V(x,y) -V(x, y+z)$, for $x\in X, y,z \in
\R$. Using (\ref{Vphi}), we obtain:
\begin{eqnarray*}
&u_z(x+\alpha, y + \varphi(x)) = V(x+\alpha, y+ \varphi(x))
- V(x+\alpha, y+z +\varphi(x)) \\
&= V(x, y )+ \varphi(x+\gamma) - [V(x, y +z )+ \varphi(x+\gamma)] =
V(x, y ) - V(x, y +z ) = u_z(x,y).
\end{eqnarray*}
Therefore $u_z$ is $\tilde T_1$-invariant, hence, by ergodicity of
$T_{\alpha, \varphi}$, for every $z$, $u_z(x, y)$ is a.e. equal to a
constant $c(z)$.

Since $u_z$ satisfies $u_{z_1+z_2} (x,y) = u_{z_1} (x,y) + u_{z_2}
(x,y+z_1)$, the previous relation implies $c(z_1+z_2) = c(z_1) +
c(z_2)$; hence, since $c$ is measurable, $c(z) = \lambda z$ for a
constant $\lambda$.

So we have for every $z$, for a.e. $(x,y)$ the relation $V(x, y+z) =
V(x, y) -\lambda z$. By Fubini it follows that for a.e. $y$, for
a.e. $(x,z)$: $V(x, y+z) = V(x, y) - \lambda z$.

Therefore, for some $y_1 \in \R$ we have $V(x, z+y_1) = V(x, y_1)
-\lambda z$; hence, setting $\psi(x) = V(x,y_1) + \lambda y_1$, we
obtain for a.e. $(x,z)$: $V(x,z) = \psi(x) - \lambda z$. Since the
Lebesgue measure is preserved on $\R$ by the map $\tilde T_2$,
necessarily $\lambda = +1$ or $\lambda= -1$.

Finally, the transformation $\tilde T_2$ has the form $(x,y) \to (x
+ \gamma, y + \psi(x))$ or  $(x,y) \to (x + \gamma, -y +
\psi(x))$.\eop

\vskip 3mm Remark that the analogous result with $\R$ replaced by
$\Z$ holds for a cocycle with values in $\Z$ which is ergodic for
the action on $X \times \Z$.

\vskip 3mm {\it Groups associated to a cocycle}

>From Equation (\ref{commute2}) it follows that $\varphi -
T_{2\gamma}{\varphi}$ is coboundary:
\begin{eqnarray}
\varphi - T_{2\gamma}{\varphi} = (\varepsilon \psi + T_{\gamma}\psi)
- T_\alpha (\varepsilon \psi + T_{\gamma} \psi). \label{commute2b}
\end{eqnarray}

Now we define several groups related to the centralizer of
$T_{\alpha, \varphi}$:
\begin{eqnarray*}
&&\Gamma:=\{\gamma: \text{ for } \varepsilon = +1 \text{ or }
\varepsilon = -1, \varepsilon \varphi - T_\gamma{\varphi} \text{ is
a coboundary for } T_\alpha \}, \\
&&\Gamma_0:=\{\gamma: \varphi - T_\gamma \varphi \text{ is a
coboundary } \psi_\gamma - T_\alpha \psi_\gamma \text{ for }
T_\alpha \}.
\end{eqnarray*}

By (\ref{commute2b}) we have $2\Gamma \subset \Gamma_0 \subset
\Gamma$. For $\gamma \in \Gamma_0$, $\psi_\gamma$ is unique up to a
constant. The family $\{\psi_\gamma, \gamma \in \Gamma_0\}$
satisfies the cocycle property on $\Gamma_0 \times X$ (up to a
constant).

For $p \in [1, \infty]$ we define
\begin{eqnarray}
\Gamma_p:=\{\gamma \in \Gamma_0 : \psi_\gamma \in L^p(\mu) \}, \ \
{\Cal C}_p (T_{\alpha, \varphi}) = \{T_{\gamma, \psi_{\gamma}},
\gamma \in \Gamma_p\}.
\end{eqnarray}

If $\gamma \in \Gamma_1$, we can choose $\psi_\gamma$ with zero
mean. The group ${\Cal C}_1 (T_{\alpha, \varphi})$ is abelian. The
cocycle property is satisfied by $\{\psi_\gamma, \gamma \in
\Gamma_1\}$: for every $\gamma, \gamma'$ in $\Gamma_1$, we have the
relation:
\begin{eqnarray*}
\psi_{\gamma' + \gamma} = \psi_\gamma + \psi_{\gamma'} (. + \gamma)
= \psi_{\gamma'} + \psi(. + \gamma').
\end{eqnarray*}

\vskip 3mm {\bf A general result on coboundaries for rotations}

Now we show that $\Gamma_0$ is a small group unless the cocycle
$\varphi$ is a coboundary, which is the degenerate case.

Let us consider the general case of rotations on a compact abelian
group $K$. For $\gamma \in K$, $T_\gamma$ denotes the rotation
(translation) by $\gamma$ on $K$. Let $T_\alpha$ be a given ergodic
rotation on $K$ defined by an element $\alpha \in K$.

The following proposition is an easy consequence of Theorem 6.2 in
\cite{MoSc80} and of the proposition p. 178 in (\cite{Le93}) (Lemma
\ref{cobQuasicob} below).

\begin{proposition}\label{cobord} Let $\varphi$ be a measurable function
on $K$. If for every $\gamma$ in a set of positive measure in $K$
there exists a measurable function $\psi_\gamma$ such that $\varphi
- T_\gamma \varphi = \psi_\gamma - T_\alpha \psi_\gamma$, then
$\varphi$ is an additive quasi-coboundary:
\begin{equation*}
\varphi = c + T_\alpha h - h,
\end{equation*}
for a measurable function $h$ and a constant $c$. If $\varphi$ is
integrable, then $c= \int \varphi \ \mu$.
\end{proposition}

\vskip 3mm
\begin{lemma} (\cite{Le93}) \label{cobQuasicob} Let
$\varphi$ be a measurable real function on $K$. If $e^{2\pi i
(\varphi - T_\gamma\varphi)}$ is a $T_\alpha$-coboundary for every
$\gamma$ in a subset of positive measure in $K$, then there are a
measurable function $\zeta_s$ of modulus 1 and $\lambda_s$ of
modulus 1 such that: $e^{2\pi i s\varphi} = \lambda_s \ T_\alpha
\zeta_s / \zeta_s$.
\end{lemma}

\vskip 3mm
\section{\bf Coboundary equations for irrational rotations \label{SectCob}}

This section is devoted to the coboundary equations over irrational
rotations, either linear equations (with Fourier's series methods)
or multiplicative equation (Guenais-Parreau's results). The
following step functions are used:

{\it Notation} \  Let $\beta$ be a fixed real number. For any real
number $\gamma$, with the notation $T_\gamma$ for the translation $x
\rightarrow x+\gamma \mod 1$, we will consider the cocycles
generated over an irrational rotation $T_\alpha$ by the step
functions
\begin{eqnarray}
\varphi_\beta = 1_{[0, \beta]} - \beta, \ \ \varphi_{\beta,
\gamma}:= 1_{[0,\beta]}-1_{[0,\beta]}(.+\gamma) = \varphi_\beta -
T_\gamma \varphi_\beta. \label{psigamma}
\end{eqnarray}

\subsection{Classical results, expansion in basis $q_n \alpha$}

\

 First of all we recall classical facts on continued fractions and
on expansion of a real $\beta$ in basis "$q_n \alpha$" (Ostrowski
expansion).

In the following $\alpha \in ]0,1[$ is an irrational number and
$[0;a_1,..., a_n,...]$ is its continued fraction expansion. Let
$(p_n/q_n)_{n \ge 0}$ be the sequence of its convergents. Recall
that $p_{-1}=1$, $p_0=0$, $q_{-1}=0$, $q_0=1$ and, for $n \ge 1$~:
\begin{equation} \label{converg_eq}
p_n = a_n p_{n-1}+p_{n-2}, \ q_n = a_n q_{n-1}+q_{n-2}, \ (-1)^n =
p_{n-1} q_n - p_n q_{n-1}.
\end{equation}

{\it Notations} \ For $u\in\R$, put $[u]$ for its integral part and
$\|u\| := \inf_{n \in \Z} |u - n|$.

For $n \ge 0$ we have $\|q_n \alpha\| = (-1)^n (q_n \alpha - p_n)$
and the following inequalities (cf. \cite{HW1}):
\begin{eqnarray}
1 &=& q_n\|q_{n+1} \alpha \| + q_{n+1} \|q_n \alpha\|, \nonumber \\
{1\over 2} {1 \over q_{n+1}} \leq {1\over q_{n+1}+q_n} &\le& \|q_n
\alpha\| \le {1\over q_{n+1}} = {1 \over a_{n+1} q_n+q_{n-1}}, \label{f_3} \\
{1\over 2} {1 \over q_{n+1}} \leq \|q_n \alpha \| &\leq& \|k
\alpha\|,\ \mbox{for}\ 1\le k < q_{n+1} \label{f_4}.
\end{eqnarray}

An irrational number $\alpha = [0;a_1,..., a_n,...]$ has {\it
bounded partial quotients} (abbreviated in "{\it is of bounded
type}") if the sequence $(a_n)$ is bounded.

\vskip 3mm {\it Expansion in basis $q_n \alpha$ (Ostrowski
expansion)}

For $\beta \in \T^1$ we consider the following representation
introduced by Ostrowski (1921)
\begin{eqnarray}
\beta = \sum_1^\infty b_j(\beta) \, q_j \,\alpha {\rm \ mod \ } 1,
\label{ostrowski}
\end{eqnarray}
where $(b_j(\beta))_{j \geq 0}$ is a sequence in $\Z$.

Any $\beta \in \T^1$ has such an expansion. If $\sum_{j \geq 1}
{|b_j(\beta)| \over a_{j+1}} < \infty$, the representation is unique
up to a finite number of terms.  It is shown in \cite{GuPa06} that
this condition is equivalent to $\sum \|q_j \beta\| < \infty$. For
$r \geq 1$, we call $H_r(\alpha)$ the subgroup
\begin{eqnarray*}
H_r(\alpha):=\bigl\{\beta = \sum_1^\infty b_j(\beta) \,q_j \,\alpha
{\rm \ mod \ } 1: \sum_{j \geq 1} {|b_j(\beta)|^r \over a_{j+1}} <
\infty\}.
\end{eqnarray*}

\vskip 3mm \goodbreak \subsection{Linear and multiplicative
equations for $\varphi_{\beta}$ and $\varphi_{\beta, \gamma}$}

\

\vskip 3mm
\subsubsection{Fourier conditions}

\

 For $\varphi(x) = \sum_{n \in \Z} \varphi_n e^{2\pi i n x}$ in $L^1(\T)$
with $\int \varphi \, d\mu = 0$, if the coboundary equation $\varphi
= h - T_\alpha h$ has a solution $h \in L^1(\T)$, the Fourier
coefficients of $h$ are $h_n = {\varphi_n \over 1- e^{2\pi i n
\alpha}}$. Therefore the necessary and sufficient condition for the
existence of a $L^2$ solution, for $\varphi \in L^2(\T)$  is
$\sum_{n \in \Z \setminus \{0\}} {|\varphi_n|^2 \over \|n
\alpha\|^2} < +\infty$.

As it is well known, under diophantine assumptions on $\alpha$ and
regularity of the function $\varphi$, the coboundary equations can
be solved. We recall briefly this fact.

The type of an irrational number $\alpha$ is $\eta \ge 1$ such that
\begin{eqnarray}
\inf_{k\neq 0} \,[k^{\eta - \varepsilon}\|k \alpha\|] =0, \ \
\inf_{k\neq 0} [ k^{\eta + \varepsilon} \|k \alpha\|] > 0, \ \forall
\varepsilon
>0. \label{type1}
\end{eqnarray}

Recall that the type of a.e. $\alpha$ is 1. From a result of V. I.
Arnold (\cite{Arn66}) (see also M. Herman (\cite{He79})), we have:
\begin{thm} (\cite{Arn66}) If $\alpha$ is of type $\eta$ and $\varphi(x)
= \sum_{n \not= 0} \varphi_n e^{2\pi i n x}$ with $\varphi_n =
O(n^{-(\eta + \delta)})$ and $\delta > 0$, then $\varphi_\gamma(x)
:= \sum_{n \not= 0} \varphi_n {1 - e^{2\pi i n \gamma} \over 1 -
e^{2\pi i n \alpha}} e^{2\pi i n x}$ is a well defined continuous
function for every $\gamma$ and the pairs $(\alpha, \varphi)$,
$(\gamma, \varphi_\gamma)$ define commuting skew products.
\end{thm}

Clearly this is a degenerate coboundary case in accordance with
Proposition \ref{cobord}, since we have a solution for every
$\gamma$. This is a motivation to consider step functions like the
function $\varphi_{\beta, \gamma}$ introduced above.

\vskip 3mm \subsubsection{The linear coboundary equation, a
sufficient condition for $\varphi_{\beta, \gamma}$} \label{lincob}

\

Now we give sufficient conditions in case of the step function
$\varphi_{\beta,\gamma}$ for the existence of a solution of the
linear coboundary equation (\ref{commute}).

Recall that the cocycle $\varphi_\beta$ is not a coboundary for
$\beta \not \in \Z\alpha + \Z$. This follows from the fact that
$e^{2\pi i \varphi_\beta} = e^{-2\pi i \beta}$, hence $e^{-2\pi i
\beta}$ is an eigenvalue of the rotation by if $\varphi_\beta$ is a
linear coboundary (cf. \cite{Pe73}). A stronger result is that the
cocycle defined by $\varphi_\beta$ over the rotation $T_\alpha$ is
ergodic if $\beta \not \in \Q \alpha + \Q$ (cf. Oren \cite{Or83}).

The Fourier coefficients of $\varphi_{\beta,\gamma} = 1_{[0, \beta]}
- 1_{[\gamma, \beta + \gamma]}$ are ${1 \over 2 \pi i n} (e^{2\pi i
n \beta} - 1) (e^{2\pi i n \gamma} - 1)$. The condition for
$\varphi_{\beta,\gamma}$ to be a coboundary with a transfer function
in $L^2(\T^1)$, i.e., such that the functional equation
$\varphi_{\beta,\gamma}= T_\alpha h - h$ has a solution $h$ in
$L^2$, is
\begin{eqnarray}
\sum_{n \not = 0} {1 \over n^2} { \|n \beta\|^2 \|n \gamma\|^2 \over
\|n \alpha\|^2} < \infty. \label{betalambda}
 \end{eqnarray}

 For the cocycle $\varphi_{\beta,
\gamma}$ the following result is proved in Appendix:
\begin{thm} \label{majBeta4} If $\beta \in \T^1$ is in $H_4(\alpha)$
then we have
$$\sum_{n \not = 0} {1 \over n^2} {\|n \beta\|^4 \over
\|n \alpha\|^2} < \infty.$$

If $\beta, \gamma$ are in $H_4(\alpha)$, then (\ref{betalambda})
holds and there is $\psi_{\beta, \gamma}$ in $L^2(\T^1)$ solution of
\begin{eqnarray}
1_{[0, \beta]} - T_\gamma1_{[0, \beta]} = \psi_{\beta, \gamma} -
T_\alpha \psi_{\beta, \gamma}. \label{solL2}
\end{eqnarray}
\end{thm}

Therefore, if $\alpha$ is not of bounded type (i.e., has unbounded
partial quotients), there is an uncountable set of pairs of real
numbers $\beta$ and $\gamma$ such that $\varphi_{\beta,\gamma}:=
1_{[0, \beta]} - 1_{[0, \beta]}(.+\gamma)$ is a coboundary $\psi -
T_\alpha \psi$ for $T_\alpha$ with $\psi$ in $L^2$.

Remark that by Shapiro's result (cf. \cite{Sh73}) on the difference
of two indicators of intervals, $\psi$ is not in $L^\infty$, unless
$\beta$ and $\gamma$ are in $\Z\alpha + \Z$.

\vskip 3mm
\subsubsection{Multiplicative equation: a necessary and
sufficient condition}

\

Now we consider the multiplicative functional equation for
$\varphi_\beta$:
\begin{equation}\label{equafonc1}
e^{2i\pi s \varphi_\beta}=e^{2i\pi t} \, {T_\alpha f \over f},
\end{equation}
$f$ is a measurable function which can be assumed of modulus 1.

Equation (\ref{equafonc1}) was studied by W. Veech \cite{Ve69}, then
by K. Merril \cite{Me85} who gave a sufficient condition on $(\beta,
s, t)$ for the existence of a solution. M. Gu\'enais and F. Parreau
have shown that this condition is sufficient and they have extended
it to more general step functions:

\begin{thm}\label{NSCdtion2} (\cite{GuPa06}, Theorems 1 and 2)
\hfill \break a) Equation (\ref{equafonc1}) has a measurable
solution $f$ for the parameters $(\beta, s, t)$ if and only if there
is a sequence of integers $(b_{n})$ such that:
\begin{eqnarray*}
&&\beta=\sum_{n\geq 0}b_{n}q_{n}\alpha \textrm{ mod }1, \text{ with
} \sum_{n\geq 0}\frac{|b_{n}|}{a_{n+1}}<\infty, \ \ \sum_{n\geq
0}\Vert b_{n}s\Vert^{2}<\infty, \\
&&t= k \alpha-\sum_{n\geq 0} [b_{n}s] \, q_{n}\alpha \textrm{ mod
}1, \text{ for an integer } k.
\end{eqnarray*}
b) Let $\varphi: \T^1 \to \R$ be a step function with integral 0 and
jumps $-s_j$ at distinct points $(\beta_j, 0\leq j\leq m$), $m\geq
1$, and let $t\in \T$. Suppose that there is a partition $\cal P$ of
$\{0,..,m\}$ such that for every $J\in {\cal P}$ and $\beta_J\in
\{\beta_j,j \in J\}$ the following conditions are satisfied:

(i) \ $\sum_{j\in J}s_j \in \Z$; \hfill \break (ii)\ for every $j
\in J$, there is a sequence of integers $(b_n^j)_n$ such that
$$\beta_j= \beta_{J}+ \sum_{n\geq 0}
b_n^jq_n\alpha \mod 1, {\ with \ } \sum_{n\geq 0}
\frac{|b_n^j|}{a_{n+1}}<+\infty, \ \sum_{n\geq 0}\Bigl\|\sum_{j \in
J}b_n^js_j\Bigr\|^2 <+\infty;$$ (iii) \ $t=k \alpha -\sum_{J\in
{\cal P}}t_J$, with $k \in \Z$ and
$$t_J=\beta_J\sum_{j\in J}s_j +\sum_{n\geq 0} \Bigl[\sum_{j \in J}
b_n^js_j \Bigr] q_n \alpha \mod 1.$$

Then there is a measurable function $f$ of modulus 1 solution of
\begin{equation}
e^{2i\pi \varphi}= e^{2i\pi t} \, T_\alpha f / f. \label{MultEquat}
\end{equation}

Conversely, when $\sum_{j \in J} s_j \notin \Z$ for every proper non
empty subset $J$ of $\{0, .., m\}$, these conditions are necessary
for the existence of a measurable solution of (\ref{MultEquat}).
\end{thm}

\vskip 3mm
\begin{rem} In the situation of Theorem \ref{majBeta4},
the multiplicative equation for $s\varphi_{\beta, \gamma}$ has a
solution for every $s \in \R$. Observe that the necessary condition
of Theorem \ref{NSCdtion2} b) does not apply to $s\varphi_{\beta,
\gamma}$ (no condition on $s$). Indeed the set of discontinuities of
$\varphi_{\beta, \gamma}$ is $J=\{0, \beta, -\gamma, \beta -
\gamma\}$ with respective jumps: $+1, -1, +1, -1$. There is a
decomposition of $J$ into $J_1=\{0, \beta\}$, $J_2=\{-\gamma, \beta
- \gamma\}$ and the sum of jumps is 0 for each of these subsets.
\end{rem}

\vskip 8mm \section{\bf Applications}

\vskip 3mm
\subsection{Non ergodic cocycles with ergodic compact quotients}
\label{nonErgCQt}

\

A first application of the results of Section \ref{SectCob} is the
construction of simple examples of non~regular cocycles with
ergodicity of all compact quotients.

By using the sufficient condition of Theorem \ref{NSCdtion2} a), we
construct non~regular (hence non ergodic) $\Z$-valued cocycles given
by the step cocycles $\varphi_{\beta, \gamma}$ defined in
(\ref{psigamma}) over rotations such that all compact quotients in
$X \times \Z/ a\Z$ are ergodic (see also \cite{Co09},
\cite{CoPi12}).

Let us recall that for every irrational number $\alpha$, for almost
every $(\beta, \gamma)$ the cocycle $\varphi_{\beta, \gamma}$ is
ergodic\footnote{See Th\'eor\`eme 5 in \cite{Co76}, where ergodicity
is proved for $T_{\alpha, \varphi}$,when $\varphi$ is a step
function, under a generic condition on the discontinuity points of
$\varphi$ called Condition (A').}. Therefore clearly we are
interested here in special, non generic, sets of values of $(\beta,
\gamma)$.

\begin{thm} \label{ergoComp} If $\alpha$ is not of bounded type,
there is $\beta$ in $H_1(\alpha)$ such that for a.e. $\gamma$:
\hfill \break a) the cocycle $\varphi_{\beta, \gamma}$ is
non~regular; \hfill \break b) all compact quotients $T_{\alpha,
\varphi_{\beta, \gamma}} \mod a: (x, y \mod a) \to (x+\alpha, y
+\varphi_{\beta, \gamma}(x) \mod a)$ are ergodic.
\end{thm}
\proof \ a) If $\alpha$ is not of bounded type, by Theorem
\ref{NSCdtion2} a) there is a non-countable set of values of $\beta$
such that, for a non-countable set of values of $s$, there are a
number $\lambda$ of modulus 1 and a measurable function $f$ of
modulus 1 such that $e^{2\pi i s \varphi_\beta} = \lambda \,
{T_\alpha f \over f}$.

We can take $\beta \not \in \alpha \Z + \Z$ and $s \not \in \Q$. For
this choice of $\beta$ and of $s$, $e^{2\pi i s(\varphi_\beta -
T_\gamma \varphi_\beta)}$ is a multiplicative coboundary for every
$\gamma$.

On the other hand, if $1_{[0,\beta]} - T_\gamma 1_{[0,\beta]}$ is an
additive coboundary for every $\gamma$ in a set of positive measure,
then by Proposition \ref{cobord} this implies that $1_{[0,\beta]} -
\beta$ is an additive coboundary which is not the case (cf.
\ref{lincob}).

Therefore for a.e. $\gamma \in \R$, $\varphi_{\beta, \gamma}$ is not
an additive coboundary. For such a value of $\gamma$, Lemma
\ref{sgroupdis} shows that ${\overline {\ev}}(\varphi_{\beta,
\gamma})= \{0, \infty \}$ and $\varphi_{\beta, \gamma}$ is
non~regular.

b) Now we construct in $H_1(\alpha)$ a more restricted set of
$\beta$ such that, for a.e. $\gamma$, the action of $T_{\alpha,
\varphi_{\beta, \gamma}}$ on the compact quotients $X \times \Z/a\Z$
are ergodic for all $a \in \Z- \{0\}$.

This done is two steps: if $\alpha$ is of non bounded type, we
construct $\beta \in H_1(\alpha)$ such that
\begin{eqnarray}
\{s : \sum_n \|b_n(\beta) s \|^2 < \infty \} \cap \Q = \Z,
\label{specIrr}
\end{eqnarray}
then show that this implies the desired property

\vskip 3mm 1) There exists a strictly increasing sequence of
integers $({j_n})$ and a sequence of integers $(d_n \geq 1)$ such
that, if one defines the subsequence $(b_{j_n})$ by
\begin{eqnarray}
b_0=1, \ b_{-1}=0, \  b_{j_{n+1}}=d_n b_{j_n}+ b_{j_{n-1}} \text{
for }  n\geq 1, \label{bnjanj}
\end{eqnarray}
then the conditions $\sum (\frac{b_{j_n}}{a_{j_{n}+1}})<\infty$ and
$\sum (\frac{b_{j_n}}{b_{j_{n+1}}})^2<\infty$ are satisfied. We
complete the sequence $(b_n)$ by zeroes.

For instance, we can choose $d_n = n$ for all $n \geq 1$ and then
$({j_n})$ such that the series $\sum { n! \over a_{j_n+1}}$
converges.

The condition $\sum (\frac{b_{j_n}}{b_{j_{n+1}}})^2<\infty$ insures
the existence of an uncountable set of values of $s$ such that
$\sum_n \|b_n s\|^2<\infty$. In particular, there is $s \not \in
\mathbb{Q}$ for which this condition holds.

Suppose that ${u \over v}$, with $u,v$ coprime integers, satisfies
$\sum_n \|b_n {u \over v} \|^2 < \infty$. For $n$ big enough, $v$
divides $u b_n$. As $b_{j_n}$ and $b_{j_{n+1}}$ are mutually coprime
(by the choice of initial values and Equation (\ref{bnjanj})), we
have $v= \pm 1$.

\vskip 3mm 2) Let $\beta$ such that $b_n(\beta) = b_n$. We have
shown above that (\ref{specIrr}) holds and the non~regularity of
$\varphi_{\beta,\gamma}$ for almost all $\gamma$. Now we prove that,
for a.e. $\gamma$, all compact quotients of $T_{\alpha,
\varphi_{\beta,\gamma}}$ are ergodic.

Let us suppose on the contrary that there is a set $D$ of positive
measure such that, for every $\gamma \in D$, there is an integer $a$
such that $T_{\alpha, \varphi_{\beta,\gamma}} \text{ mod } a$ is non
ergodic.

Using Fourier series representation of $T_{\alpha, \varphi_{\beta,
\gamma}}$-invariant $a$-periodic functions, this would imply the
following: there are integers $a$ and $k$, with $a, k$ coprime, and
a set $D_{a,k}$ of positive measure such that for every $\gamma \in
D_{a,k}$ there exists a measurable function $f_\gamma$ satisfying:
\begin{equation}
e^{-2i\pi \frac ka (\varphi_\beta - T_\gamma \varphi_\beta) (x)} =
f_\gamma(x)/f_\gamma(x+\alpha). \label{equak}
\end{equation}

Lemma \ref{cobQuasicob} implies the existence of $t$ and $h$ such
that
$$e^{-2i\pi \frac ka \varphi_{\beta}(x)} = e^{2i\pi t} \, h(x)/h(x+\alpha).$$
As the conditions in Theorem \ref{NSCdtion2} a) are necessary, this
implies that $\sum \|b_n \frac ka\|^2<\infty$, contrary to
(\ref{specIrr}).

Remark that, by strengthening the conditions in the construction of
$\beta$, we can also find $\beta \in H_4(\alpha)$ with the previous
properties. For such a $\beta$, by Theorem \ref{majBeta4} there is
an uncountable set of values of $\gamma$ for which
$\varphi_{\beta,\gamma}$ is a coboundary. \eop

\vskip 5mm \goodbreak \subsection{Examples of non trivial and
trivial centralizer}

\

The results of Subsection \ref{thmCommut} lead to the following
questions for a given rotation $T_\alpha$ and a function $\varphi$:
\hfill \break - for which $\gamma \in \T^1$ is there a solution to
the commutation equation $\varphi - T_\gamma \varphi = \psi -
T_\alpha \psi$? \hfill \break - what is the centralizer of
$T_{\alpha, \varphi}$?

In this subsection, from Theorem \ref{majBeta4} we obtain that the
centralizer of $T_{\alpha, \varphi_\beta}$ is non countable for
$\beta \in H_4(\alpha)$. Then we show that the centralizer ${\Cal
C}(T_{\alpha, \varphi_\beta})$ is also non trivial when $\beta \in
H_1(\alpha)$. In a second part, we investigate a property of
``rigidity'' for $\alpha$ of bounded type, with an example of a
small centralizer.

\vskip 3mm \subsubsection{ Case of a non trivial centralizer}

\

Let $\alpha$ be an irrational number which is not of bounded type
and $\beta$ a real number. Let us consider $\varphi = \varphi_\beta
= 1_{[0, \beta]} - \beta$.

If $T_{\alpha, \varphi_\beta}$ is ergodic, by Theorem
\ref{thmCommut} and the commutation relation (\ref{commute2b}), the
square of the elements of ${\Cal C}(T_{\alpha, \varphi_\beta})$ are
of the form $T_{\gamma, \psi}$ with $\psi$ a measurable function and
$\gamma$ such that $1_{[0, \beta]}(.) - 1_{[0, \beta]} (.+2\gamma) =
\psi - T_\alpha \psi$.

By Theorem \ref{majBeta4}, if $\beta$ is in $H_4(\alpha)$, the group
$\Gamma_2$ defined in Subsection \ref{subsect-centralz} contains the
group $H_4(\alpha)$, which is a non~countable group if $\alpha$ is
not of bounded type.

Now we would like to weaken the condition on $\beta$ and still get a
non trivial centralizer. It is interesting to investigate the
properties of the cocycle $\varphi_{\beta, \beta}$ or more generally
$\varphi = a1_{[0, \beta]} - 1_{[0, a\beta]}$ with $a$ a positive
integer. This is a special situation where one can conclude that the
cocycle is a coboundary by using the result of Gu\'enais and Parreau
mentioned above.

\begin{prop} \label{exCob} If $a$ is a positive integer, the cocycle $\varphi =
a1_{[0, \beta]} - 1_{[0, a\beta]}$ is a coboundary if and only if
$\beta$ is in $H_1(\alpha)$.
\end{prop} \proof \ With the notation of Theorem \ref{NSCdtion2},
the discontinuities of $\varphi = a1_{[0, \beta]} - 1_{[0, a\beta]}$
are at $\beta_0 = 0, \beta_1 = \beta, \beta_2 = \gamma = a \beta$,
with jumps respectively $a-1, -a, 1$, we have $m= 2$ and the
partition ${\cal P}$ is the trivial partition with the single atom
$J = \{0,1,2\}$. We have $\beta_J = 0$, $\sum_{j\in J} s_j = 0$.

Suppose that $\beta \in H_1(\alpha)$ with an expansion in basis
$(q_n \alpha)$ given by
\begin{eqnarray}
\beta= \sum_{n\geq 0} b_n q_n\alpha \mod 1, {\rm \ with \ }
\sum_{n\geq 0} {|b_n| \over a_{n+1}}<+\infty, \ b_n \in \Z.
\label{betaExpan}
\end{eqnarray}
We can take $b_n^0 = 0, b_n^1 = b_n, b_n^2 = ab_n$, so that $\sum_{j
\in J} b_n^j s_j = ab_n - ab_n = 0$. For every real $s$ the
multiplicative equation $e^{2\pi i s \varphi} = T_\alpha f/f$ has a
solution. By using Theorem 6.2 in \cite{MoSc80}, we conclude that
$\varphi$ is a measurable coboundary (another proof based on the
tightness of the cocycle (that is, the tightness of the family
$(\varphi_n, n \geq 0)$) can also be given).

Conversely, if $\varphi$ is a measurable coboundary, then $e^{2\pi i
s \varphi} = T_\alpha f/f$ has a solution for every real $s$, and
this implies that $\beta$ has an expansion like in (\ref{betaExpan})
(Theorem \ref{NSCdtion2} b), necessary condition). \eop

\vskip 3mm Under the assumption $\beta \in H_1(\alpha)$ which is
weaker than the assumption of Theorem \ref{majBeta4}, Proposition
\ref{exCob} implies:
\begin{cor} \label{Cor2beta} If $\beta \in H_1(\alpha)$, the centralizer
${\Cal C}(T_{\alpha, \varphi_\beta})$ contains a non trivial element
$T_{\beta, \psi_\beta}$, where $\psi_\beta$ is a measurable function
solution of $1_{[0, \beta]} - 1_{[\beta, 2\beta]} = \psi_\beta -
T_\alpha \psi_\beta$.
\end{cor}

\vskip 3mm \begin{rem} {\rm  We have seen in the previous
considerations that, under some assumption on the expansion of
$\beta$ in basis $q_n \alpha$, the cocycle $\varphi_{\beta, \beta} =
1_{[0, \beta]} - T_\beta 1_{[0, \beta]}$ is a coboundary for the
rotation by $\alpha$, with a transfer function in a certain space:}
\end{rem}
({\it i}) \, if $\varphi_{\beta, \beta}$ is a coboundary in the
space of bounded functions, then $\beta \in \Z\alpha + \Z$ (cf.
Shapiro's result);

({\it ii}) \, if $\sum b_k^4/a_{k+1} < \infty$, then
$\varphi_{\beta, \beta}$ is a coboundary with a transfer function in
$L^2$ (see Theorem \ref{majBeta4});

({\it iii}) if $\sum |b_k|/a_{k+1} < \infty$, then $\varphi_{\beta,
\beta}$ is a coboundary with a measurable transfer function.
(Proposition \ref{exCob}). This is also necessary by Theorem
\ref{NSCdtion2} b).

\vskip 3mm \subsubsection{Example of trivial centralizer}

\

Now, for $\alpha$ of bounded type, we show the triviality of the
centralizer in the special case $\beta = \frac12$.

\begin{thm} Let $\alpha$ be of bounded type. For $\beta = \frac12$,
the centralizer of $T_{\alpha, \varphi_{\beta }}$ (acting on $X
\times \frac12 \Z$) reduces to the translations on the fibers $(x,y)
\to (x, y + \lambda)$, for a constant $\lambda \in \R$, the map
$(x,y) \to (x + \frac12, -y)$ and the powers of $T_{\alpha,
\varphi_{\beta}}$. \end{thm} \proof \ The cocycle $\varphi=
\varphi_{\frac12, \frac12} = 2\varphi_\frac12$ is known to be
ergodic as a cocycle with values in $\Z$, for every irrational
rotation (\cite{CoKe76}).

According to Theorem \ref{thmCommut} and the commutation relation
(\ref{commute2}), we consider the cocycle $u_\gamma:= \varepsilon
\varphi - T_\gamma \varphi$, where $\varepsilon$ is the constant
$+1$ or $-1$. Suppose that $\alpha$ is of bounded type and $\gamma
\not \in \Z \alpha + \Z$.

Assume that $\gamma \not =\frac12 \mod 1$, so that $u_\gamma$ has
effective discontinuities for $x = 0, \frac12, - \gamma, \frac12 -
\gamma$.

By Lemma 2.3 and Theorem 3.8 in \cite{CoPi12} to which we refer for
more details, the cocycle $u_\gamma$ satisfies a property of
separation of its discontinuities along a subsequence of
denominators of $\alpha$ and therefore its discontinuities belong to
the group of its finite essential values.

This implies that $u_\gamma$ has a non trivial essential value,
hence is not a coboundary.

The case $\gamma = \frac12 \mod 1$ corresponds to the special map
$(x,y) \to (x + \frac12, -y)$ which yields an element in the
centralizer due to the relation satisfied here: $-\varphi(x) =
\varphi(x + \frac12)$.

It remains to examine the case $\gamma = p \alpha$ $\mod 1$, with $p
\not = 0$ in $\Z$.

Suppose that $\varepsilon = -1$. Then $\varphi + T_{p\alpha}\varphi$
is a $T_\alpha$-coboundary, hence also $\varphi$, since $\varphi -
T_{p\alpha}\varphi = (\varphi + ... + T_{(p-1)\alpha}\varphi) -
T_\alpha(\varphi + ... + T_{(p-1)\alpha} \varphi)$ is a coboundary.

Since $\varphi$ is not a coboundary, necessarily $\varepsilon = +1$.

For $\varepsilon = +1$ and $\gamma = p \alpha + \ell$, we find the
powers of the map $T_{\alpha, \varphi_{\beta }}$. \eop

\vskip 3mm
\subsection{Example of a non trivial conjugacy in a group family} \label{counterex}

\

Another application is a conjugacy problem for a family of closed
subgroup over a dynamical system.

We consider the following data: a dynamical system $(X, \mu ,T)$, a
measurable family $(H_x)_{x \in X}$ of closed subgroups of a (non
commutative) topological group $G$ and a measurable function $\Phi:
X \to G$ such that the following conjugacy equation holds:
\begin{equation}
H_{T x} = \Phi(x) \ H_{x} \ (\Phi(x))^{-1}, {\rm \ for \ }
\mu_\chi{\rm-a.e. \ } x \in X. \label{conjug}
\end{equation}

We would like to give a simple example of construction of such a
family which is {\it not conjugate to a fixed closed subgroup} of
$G$ (cf. \cite{CoRa09}), i.e., such that there is no subgroup $H
\subset G$ and no measurable function $\zeta: X \to G$ solution of
the equation
\begin{equation}
H_x = \zeta(x)^{-1} H \zeta(x). \label{conjugH}
\end{equation}

Let $\theta$ be a fixed irrational number and let $G$ be the
solvable group obtained as the semi-direct product of $\mathbb{R}$
and $\mathbb{C}^2$, with the composition law:
$$(t,z_1,z_2) * (t',z_1',z_2') = (t+t',z_1+e^{2\pi i t}
z_1',z_2+e^{2\pi \theta i t}z_2').$$

The conjugate of $(0,z_1,z_2)$ by $a = (s,v_1,v_2)$ in $G$ is:
\begin{equation} (s,v_1,v_2) (0,z_1,z_2) (s,v_1,v_2)^{-1} =
(0,e^{2\pi i s} z_1,e^{2\pi \theta i s}z_2). \label{conjugue}
\end{equation}

Consider the dynamical system defined by an irrational rotation $T :
x \rightarrow x+ \alpha \mbox{ mod }1$ on $X = \T^1$. Let $\Phi : X
\rightarrow G$ be the cocycle defined by $\Phi(x) = (\varphi(x), 0,
0)$, where $\varphi$ has its values in $\Z$.

\vskip 3mm Let $H_x :=\{(0, vz_1, ve^{2\pi i \psi(x)} z_2), v \in
\mathbb{R} \}$, where $\psi$ is a measurable real function defined
below and $z_1$, $z_2$ are given real numbers. For every $x \in X$,
$H_x$ is a closed subgroup of $G$. Let us consider the function $x
\rightarrow H_x$ with values in the set of closed subgroups of $G$.
It satisfies the conjugacy relation (\ref{conjug}) if and only if
$\varphi$ has integral values and satisfies
\begin{equation}
\theta \ \varphi(x) + \psi(x) = \psi(T x) \ \hbox{ mod } 1.
\label{phipsi}
\end{equation}

Let us take $\varphi = \varphi_{\beta, \gamma} = 1_{[0,\beta]} -
1_{[0,\beta]}(. + \gamma)$. We have seen that, for every $\alpha$
which is not of bounded type, there are real numbers $\beta$ and
$\gamma$ for which the function $\varphi_{\beta, \gamma}$ is not a
coboundary and $e^{2\pi i \theta \varphi_{\beta, \gamma}}$ is a
multiplicative coboundary for some irrational values of $\theta$.

It means that for these values of the parameters, there is $\psi$
such that (\ref{phipsi}) is satisfied

\begin{prop} For these choices of $\beta, \theta$, $\varphi
= \varphi_{\beta, \gamma}$ and $\psi$, there is no subgroup $H$ such
that the equation (\ref{conjugH}) has a measurable solution $\zeta$.
\end{prop} \Proof \ \ Suppose that there are a fixed subgroup $H$ and a
measurable function $\zeta : X \rightarrow G$ solution of
(\ref{conjugH}). According to (\ref{conjugue}), this is equivalent
to the existence of a function $\rho$ defined on $X$ such that the
set
$$\{(0, ve^{2\pi i \rho(x)} z_1, ve^{2\pi i (\theta \rho(x) + \psi(x))} z_2), v
\in \mathbb{R} \}$$ does not depend on $x$. This implies that $\rho$
and $\psi + \theta \rho$ have a fixed value mod 1. Therefore
$\rho(x) - \rho(Tx) \in \Z$, $\theta (\varphi(x) - \rho(x) + \rho(T
x)) = \theta \varphi(x) + \psi(x) - \psi(T x)$ and according to
(\ref{phipsi}) the common value $\hbox{ mod } 1 $ is 0.

As $\varphi$ has integral values and $\theta$ is irrational, it
follows that $\varphi = T \rho - \rho$, contrary to the fact that
$\varphi$ is not a coboundary. \eop

\vskip 3mm \vskip 3mm \section{\bf Appendix: proof of Theorem
\ref{majBeta4}\label{append}}

For the proof of Theorem \ref{majBeta4} we need some preliminary
results. In what follows, $C$ will denote a generic constant which
may change from a line to the other.

\vskip 3mm {\it Bounds for $\|q_n \beta\|$}

Let $\beta \in [0, 1]$ be such that
\begin{eqnarray} \beta =
\sum_1^\infty b_i q_i \alpha \moda, {\rm \ with \ }\sum_1^\infty
{|b_i| \over a_{i+1}} = C_1 < \infty. \label{defbeta}
\end{eqnarray}

In the following computations, we assume that there is infinitely
many $i$'s with $b_i \not = 0$. We can assume $b_i \geq 0$.

The quantities $\|q_n \beta\|$ and $|b_n|/a_{n+1}$ are of the same
order. For all $r \geq 1$ such that $b_r \not = 0$, the following
upper bounds hold:
\begin{eqnarray*}
&&\sum_{j=1}^r b_j q_j \leq q_r(b_r + {b_{r-1} \over a_{r}} +
{b_{r-2} \over a_{r} a_{r-1}}+ ...+ {b_{1} \over a_{r}
a_{r-1}...a_2}) \leq q_r (b_r + C_1) \leq (C_1+1) b_r q_r, \\
&&\sum_{j=r}^\infty b_j \|q_j \alpha\| \leq {b_r \over q_{r+1}} +
{b_{r+1} \over q_{r+2}}+ ... \leq {1 \over q_{r+1}} (b_r + {b_{r+1}
\over a_{r+2}} + {b_{r+2} \over a_{r+2} a_{r+3}} + ...) \\
&&\quad \quad \quad \quad \quad \quad \leq {1 \over q_{r+1}} (b_r +
C_1) \leq (C_1+1) {b_r \over q_{r+1}}.
\end{eqnarray*}

\vskip 3mm For $n \geq 1$, let $\ell(n)$ be the greatest index $i
\leq n-1$ such that $b_i \not = 0$, and $m(n)$ the smallest index $i
\geq n$ such that $b_i \not = 0$. For all $r, k \geq 1$, we have
$$\|k \beta \| = \|\sum_1^\infty b_i q_i k \alpha\|
\leq \min (1,\, 2 \max (\|k \alpha \| \,\sum_1^{r-1} b_i q_i , \, k
\,\sum_{r}^\infty b_i \|q_i \alpha\|));$$ hence with $C = 2(C_1+1)$:
\begin{eqnarray}
\|k \beta \| \leq \min (1,\, C \max (\|k \alpha \| \,b_{\ell(n)}
q_{\ell(n)},\, k {b_{m(n)} \over q_{m(n)+1}})), \ \forall n,k \geq
1. \label{majkbeta}
\end{eqnarray}

Observe that since $b_{\ell(n)}$ is non zero integer,
\begin{eqnarray}
\sum_n {1 \over a_{\ell(n)+1}} \leq \sum_n {|b_{\ell(n)}| \over
a_{\ell(n)+1}} < C_1. \label{bln1}
\end{eqnarray}

We will use also that if $s$ is an integer $\geq 1$, then
\begin{eqnarray*}
&&\sum_j {b_j^{s} \over a_{j+1}} < +\infty \Rightarrow \sum_n
{b_{\ell(n)}^{s} q_{\ell(n)} \over q_{\ell(n)+1}} < +\infty.
\end{eqnarray*}

\vskip 3mm {\it Denjoy-Koksma inequality} (cf. \cite{He79})

We denote by $V(f)$ the variation of a BV (bounded variation)
function $f$ on $X = \R / \Z$, for instance a step function with a
finite number of discontinuities. If $p/q$ is a irreducible fraction
such that $\|\alpha - p/q\| < {1 / q^2}$, then for every $x\in X$
the following inequality holds:
\begin{eqnarray}
|\sum_{\ell = 0}^{q-1}f(x+\ell \alpha) - q \int f \, dy| \le V(f).
\label{f_8}
\end{eqnarray}

Let $S_n f = \sum_{k=0}^{n-1} T_\alpha^k f$ be the Birkhoff sums of
$f$ for the rotation $T_\alpha$. Using Inequality (\ref{f_8})
implies for the denominators $q_n$ of $\alpha$:
\begin{eqnarray}
\|S_{q_n}f \|_\infty \leq |\mu(f)|\, q_n + V(f), \forall n \in \N.
\label{majgamma}
\end{eqnarray}

\vskip 3mm
\begin{lem} If $f$ is a nonnegative BV function, we have,
\begin{eqnarray}
\sum_{k=q_n}^\infty {f(k\alpha) \over k^2} \leq 2 ( {\mu(f) \over
q_n} +  {V(f)  \over q_n^2}), \, \forall n \geq 1. \label{maj1}
\end{eqnarray}
\end{lem} \proof \ The inequality (\ref{majgamma}) implies
\begin{eqnarray*} \sum_{k=q_n}^\infty {f(k\alpha) \over k^2} &&\leq
\sum_{j=1}^\infty {1\over (j q_n)^2} \sum_{p=0}^{q_n-1} f((j q_n +p) \alpha) \\
&&\leq {1\over q_n^2} (\sum_{j=1}^\infty {1\over j^2})\, (\mu(f)\,
q_n + V(f)) = 2 ({\mu(f) \over q_n} + { V(f) \over
q_n^2}).\end{eqnarray*} \eop

\vskip 3mm For all $p \geq 1$, by (\ref{maj1}) applied with $f(x) =
{1 \over x^2} 1_{[{1\over p}, {1\over 2}]}(|x|)$, then applied with
$f(x) = 1_{[-{1\over p}, {1\over p}]}(x)$, we get
\begin{eqnarray}
&&\sum_{\{k \geq q_n, \, \|k \alpha \| \geq 1/p\}} {1 \over k^2} {1
\over \|k \alpha \|^2} \leq C ({ p \over q_n} + { p^2 \over
q_n^2}) \label{maj2}, \\
&&\sum_{\{k \geq q_n, \, \|k \alpha \| \leq 1/p\}} {1 \over k^2}
\leq C({1 \over q_n p} + {1 \over q_n^2}). \label{maj3}
\end{eqnarray}

\vskip 3mm On the other hand, we have, from (\ref{f_8}):
\begin{eqnarray*}
\sum_{\{0 < k < q_n, \, \|k \alpha \| \geq 1/p\}} {1 \over \|k
\alpha \|^2} && = \sum_{\{0 < k < q_n \}} {1 \over \|k
\alpha \|^2} \, 1_{[1/p, \, 1 - 1/p]} (\{k \alpha\}) \\
&&\leq 2 p^2 \sum_{\ell=1}^{[{p+1 \over 2}]} {1 \over \ell^2}
\sum_{\{0 < k < q_n \}} 1_{[{\ell \over p}, {\ell+1 \over p}[} (\{k
\alpha\})\ \leq 2 p^2 ({q_n \over p} +2) \, \sum_{\ell \geq 1} {1
\over \ell^2};
\end{eqnarray*}
hence:
\begin{eqnarray} \sum_{\{0 < k
< q_n, \, \|k \alpha \| \geq 1/p\}} {1 \over \|k \alpha \|^2} \leq C
p (q_n +p) \label{maj4}.
\end{eqnarray}

\vskip 2mm
\goodbreak
\begin{lem} \label{1-4} a) There is a finite constant $C$ such that,
for every $n \geq 1$,
\begin{eqnarray}
\sum_{k=1}^{q_n-1} {1 \over \|k \alpha \|^2} \leq C q_n^2.
\label{sumInvSqr}
\end{eqnarray}

b) For all $s \geq 1$, there exists at most one value of $k$ of the
form $k= s q_n + r$, with $r \in [1, q_n[$ such that $\| k \alpha \|
< {1 \over 4} {1 \over q_n}$, and this value satisfies $k \geq {1
\over 4} q_{n+1}$.
\end{lem} \Proof \ \ a) If $\alpha > {p_n \over q_n}$, then each interval
$[{j \over q_n}, {j+1 \over q_n})$, $1 \leq j \leq q_n-1$ contains
exactly one number of the form $\{k \alpha\}$, with $1 \leq k \leq
q_n-1$. Therefore we have
\begin{eqnarray}
\sum_{k=1}^{q_n-1} {1 \over \|k \alpha \|^2} \leq \sum_{\ell =
1}^{q_n-1} (\ell / q_n)^{-2} \leq q_n^2 \sum_{\ell =1}^{+\infty}
\ell^{-2}. \label{sumInvSqr1}
\end{eqnarray}
When $\alpha < {p_n \over q_n}$, the same is true for $j=1, \ldots,
q_n-2$. Furthermore there is an exceptional value $k_1$ (the value
such that $k_1 p_n = 1 \text{ mod } q_n$) for which $0 < \{k_1
\alpha\} < {1 \over q_n}$. By (\ref{f_4}) we know that $\|k \alpha
\| \geq {1 \over 2 q_n}$ for $1 \leq |k| < q_n$. Therefore ${1 \over
2 q_n} < \{k_1 \alpha\} < {1 \over q_n}$ which add a contribution of
$4 q_n^2$ in (\ref{sumInvSqr1}). This implies (\ref{sumInvSqr1}).

b) For a given integer $s \geq 1$, suppose that there are two
different values of the form $k_i = s q_n + r_i$, with $r_i \in [1,
q_k[$, $i=1,2$, satisfying: $\| k_i \alpha \| < {1 \over 4} {1 \over
q_n}$.

Then we have, for some $r_0 \in [1, q_n[$, $\|r_0 \alpha \| < {1
\over 2} {1 \over q_n}$, which contradicts that for $r_0 \in [1,
q_n[$ we have $\| r_0 \alpha \| \geq \| q_{n-1} \alpha \| \geq
{1\over 2 q_n}$ by (\ref{f_3}) and (\ref{f_4}).

\vskip 3mm Let $k= s q_n + r$, with $r \in [1, q_n[$ and $\| k
\alpha \| < {1 \over 4} {1 \over q_n}$. Put $\lambda = s q_n /
q_{n+1}$. The condition $\| k \alpha \| < {1 \over 4} {1 \over q_n}$
implies
$$\| r \alpha \| < {1 \over 4 q_n} + {s \over q_{n+1}} \leq (\lambda
+ {1\over 4}) {1 \over q_n}.$$

\vskip 3mm As $\| r \alpha \| \geq \|q_{n-1} \alpha \| \geq {1 \over
2 q_n}$, for $r \in [1, q_n[$, we get $\lambda > {1 \over 4}$. \eop

\vskip 3mm The proof of Theorem \ref{majBeta4} relies on the
expansion of $\beta$ in basis $q_n \alpha \text{ mod } 1$. We
suppose that $\beta \in [0,1[$ satisfies (\ref{defbeta}), so that we
can apply (\ref{majkbeta}).

We denote by $J$ and $J'$ the sets of integers defined by
\begin{eqnarray}
&&J:= \{k = sq_n, s= 1, ..., a_{n+1}, n=1, 2, ... \}, \\
&&J':= \bigcup_{n=1}^\infty ([q_n, q_{n+1}[ \cap \{k: \|k \alpha\| <
{1 \over 4 q_n } \}).
\end{eqnarray}

\vskip 3mm
\begin{lem} \label{hypoBeta2}
\begin{eqnarray} \sum {b_j(\beta)^2 \over a_{j+1}} < \infty
\ \Rightarrow \sum_{n \not = 0, n \not \in J} {1 \over n^2} {\|n
\beta\|^2 \over \|n \alpha\|^2} < \infty.
\end{eqnarray}
\end{lem} \proof \ Up to a constant factor, we have
\begin{eqnarray*}
&&\sum_{k \not = 0, k \not \in J} {\|k \beta\|^2 \over k^2 \|k
\alpha\|^2} = \sum_{n=0}^\infty \ \sum_{k \not \in J, \, q_n \leq k
< q_{n+1}}
{\|k \beta\|^2 \over k^2 \|k \alpha\|^2}\\
&& \leq \sum_{n=0}^\infty \ \sum_{q_n \leq k < q_{n+1}, \, \|k
\alpha\| \geq 1/q_{\ell(n)}} {\|k \beta\|^2 \over k^2 \|k
\alpha\|^2} + \sum_{n=0}^\infty \ \sum_{k \not \in J, \, q_n \leq k
< q_{n+1}, \, \|k \alpha\| <
1/q_{\ell(n)}} {\|k \beta\|^2 \over k^2 \|k \alpha\|^2} \\
&& \leq \sum_{n=0}^\infty \ \sum_{q_n \leq k < q_{n+1}, \, \|k
\alpha\| \geq 1/q_{\ell(n)}} {1 \over k^2 \|k \alpha\|^2} +
\sum_{n=0}^\infty \ \sum_{k \not \in J, \, q_n \leq k < q_{n+1}, \,
\|k \alpha\| <
1/q_{\ell(n)}} {\|k \beta\|^2 \over k^2 \|k \alpha\|^2} \\
&&\leq (A) + (B) + (C) + (D),
\end{eqnarray*}
with (using (\ref{majkbeta}) for (B) and (C)):
\begin{eqnarray*}
&&(A) := \sum_n \sum_{q_n \leq k < q_{n+1}, \, \|k \alpha\| \geq
1/q_{\ell(n)}} {1\over k^2 \|k \alpha\|^2},\\
&&(B) := \sum_n \sum_{q_n \leq k < q_{n+1}, \, \|k \alpha\| <
1/q_{\ell(n)}} b_{\ell(n)}^2 q_{\ell(n)}^2 {1\over
k^2}, \\
&&(C) := \sum_n \sum_{{k \not \in J\cup J', \, q_n \leq k <
q_{n+1}}} b_{m(n)}^2 {1 \over q_{m(n)+1}^2} {1\over \|k \alpha\|^2}. \\
&&(D) := \sum_{{k \not \in J}, \, {k \in J'}} {\|k \beta\|^2\over
k^2 \|k \alpha\|^2}.
\end{eqnarray*}

Observe that $q_n \geq q_{\ell(n)+1} > a_{\ell(n)+1} \,
q_{\ell(n)}$, since $\ell(n)+1 \leq n$. We have from (\ref{maj2})
and from (\ref{maj3}) applied with $p = q_{\ell(n)}$ and from
(\ref{bln1}):
\begin{eqnarray*}
(A) &&\leq C \sum_n ({q_{\ell(n)} \over q_n} + {q_{\ell(n)}^2 \over
q_n^2}) \leq C\sum_j ({1 \over a_{\ell(n)+1}} + {1 \over
a_{\ell(n)+1}^2})
< \infty,\\
(B) &&\leq C \sum_n b_{\ell(n)}^2 q_{\ell(n)}^2
({1\over q_n q_{\ell(n)}} + {1 \over q_n^2}) \\
&&\leq C \sum_n b_{\ell(n)}^2 ( {q_{\ell(n)} \over q_n} +
{q_{\ell(n)}^2 \over q_n^2}) \leq \sum_j {b_j ^2 \over a_{j+1}} (1 +
{1 \over a_{j+1}}) < +\infty,
\end{eqnarray*}
and from (\ref{maj4}), as the sum is taken over indices $k \not \in
J'$, i.e. such that $\|k \alpha \| \geq {1\over 4 q_n}$~:
\begin{eqnarray*}
(C) &&\leq C \sum_n b_{m(n)}^2 {4 q_n(q_{n+1} + 4 q_n) \over
q_{m(n)+1}^2} \\
&&\leq C \sum_n b_{m(n)}^2 {q_n \over q_{m(n)+1}} \leq C' \sum_n
{b_{m(n)}^2 \over a_{m(n)+1}} \leq C' \sum_j {b_j^2 \over a_{j+1}} <
+\infty.
\end{eqnarray*}

We are left with the convergence of the series (D). By
(\ref{majkbeta}) it suffices to prove the following convergence
\begin{eqnarray*}
&&(E) := \sum_n \sum_{k \not \in J, k \in J', \, q_n \leq k <
q_{n+1}} b_{\ell(n+1)}^2 q_{\ell(n+1)}^2 {1 \over k^2} < \infty,\\
&& (F) := \sum_n \sum_{k \not \in J, k \in J', \, q_n \leq k <
q_{n+1}} b_{m(n+1)}^2 {1 \over q_{m(n+1)+1}^2} {1 \over \|k \alpha
\|^2} < \infty.
\end{eqnarray*}
By Lemma \ref{1-4}, we have:
\begin{eqnarray*}
(E) \leq C \sum_n b_{\ell(n+1)}^2 \, q_{\ell(n+1)}^2 {a_{n+1} \over
q_{n+1}^2} \leq C \sum_n b_{\ell(n+1)}^2 {q_{\ell(n+1)} \over
q_{\ell(n+1)+1}} < \infty.
\end{eqnarray*}

To bound (F), we use Lemma \ref{1-4} a):
\begin{eqnarray*}
(F) && \leq C \sum_n b_{m(n+1)}^2 {1 \over q_{m(n+1)+1}^2} q_{n+1}^2
(\sum_{\ell = 1}^\infty {1 \over \ell^2}) \\
&&\leq C \sum_n b_{m(n+1)}^2 {q_{m(n+1)}^2 \over q_{m(n+1)+1}^2}
(\sum_{\ell = 1}^\infty {1 \over \ell^2}) \leq C \sum {b_j^2 \over
a_{j+1}^2} <\infty.
\end{eqnarray*}
\eop

\vskip 3mm {\bf Proof of Theorem \ref{majBeta4}} Let $\beta \in
H_r(\alpha)$, i.e.,
\begin{eqnarray}
\sum {b_j(\beta)^4 \over a_{j+1}} < \infty, \label{hypoBeta4}
\end{eqnarray}
Taking into account Lemma \ref{hypoBeta2}, it remains to show the
convergence of
$$\sum_{n \not = 0} \sum_{s= 1}^{a_{n+1}} {\|sq_n \beta\|^4
\over s^2 q_n^2 \|sq_n \alpha\|^2}.$$

By (\ref{majkbeta}) applied with $k = q_n$, it suffices to prove the
convergence of the series
\begin{eqnarray*}
&&(G) := \sum_n b_{\ell(n)}^2 q_{\ell(n)}^2 {1 \over q_n^2}
(\sum_{s=1}^{a_{n+1}} {1 \over s^2}), \\
&&(H) := \sum_n b_{m(n)}^4 {1 \over q_{m(n)+1}^4} (\sum_{s=
1}^{a_{n+1}} q_n^2 q_{n+1}^2).
\end{eqnarray*}

Since $n \leq m(n)$ and $\ell(n) + 1 \leq n$, we have from
(\ref{hypoBeta4}):
\begin{eqnarray*}
(G) &&\leq C \sum_n b_{\ell(n)}^2 {q_{\ell(n)}^2 \over q_n^2} \leq C
\sum_n b_{\ell(n)}^2 {q_{\ell(n)}^2 \over q_{\ell(n)+1}^2} \leq C
\sum_n {b_{\ell(n)}^2 \over a_{\ell(n)+1}^2}
\leq \sum_j {b_{j}^2 \over a_{j+1}^2} \leq \sum_j {b_{j}^2 \over a_{j+1}} < \infty, \\
(H) &&\leq C \sum_n b_{m(n)}^4 {a_{n+1} q_n^2 q_{n+1}^2 \over
q_{m(n)+1}^4} \leq C \sum_n b_{m(n)}^4 {q_n q_{n+1}^3 \over
q_{m(n)+1}^4} \leq C \sum_n b_{m(n)}^4 {q_{m(n)} q_{m(n)+1}^3 \over
q_{m(n)+1}^4} \\
&&= C \sum_n b_{m(n)}^4 {q_m(n) \over q_{m(n)+1}} \leq C \sum_n
{b_{m(n)}^4 \over a_{m(n)+1}} \leq \sum_j {b_{j}^4 \over a_{j+1}} <
\infty.
\end{eqnarray*}

Therefore, if $\beta \in H_4(\alpha)$, we have $\sum_{n \not = 0} {1
\over n^2} {\|n \beta\|^4 \over \|n \alpha\|^2} < \infty$.

 For the second statement of Theorem \ref{majBeta4}, observe that, if
$\beta$ and $\gamma$ belong to $H_4(\alpha)$, by the previous
inequality and Cauchy-Schwarz inequality
\begin{eqnarray}
\sum_{n \not = 0} {1 \over n^2} { \|n \beta\|^2 \|n \gamma\|^2 \over
\|n \alpha\|^2} < \infty. \label{CSbetagamma}
\end{eqnarray}

Recall that the Fourier coefficients of $\varphi_{\beta,\gamma}$ are
${1 \over 2 \pi i n} (e^{2\pi i n \beta} - 1) (e^{2\pi i n \gamma} -
1)$. The condition for $\varphi_{\beta,\gamma}$ to be a coboundary
with a transfer function in $L^2(\T^1)$, i.e., such that the
functional equation $\varphi_{\beta,\gamma}= T_\alpha h - h$ has a
solution $h$ in $L^2$, is fulfilled by (\ref{CSbetagamma}). \eop

IRMAR, CNRS UMR 6625, Universit\'e de Rennes I,\\
Campus de Beaulieu, 35042 Rennes Cedex, France

conze@univ-rennes1.fr \\
jonat.marco@gmail.com
\end{document}